\documentclass[11pt]{article}
\usepackage{latexsym,amsmath,amsfonts,hyper,setspace}
\def\nbw{section}

\numberwithin{equation}{\nbw}

\newtheorem{theorem}{Theorem}
\numberwithin{theorem}{\nbw}
\newtheorem{proposition}{Proposition}
\numberwithin{proposition}{\nbw}
\newtheorem{corollary}{Corollary}
\numberwithin{corollary}{\nbw}

\numberwithin{lemma}{\nbw}

\numberwithin{question}{\nbw}

\numberwithin{conjecture}{\nbw}

\numberwithin{assumption}{\nbw}

\numberwithin{definition}{\nbw}

\numberwithin{notation}{\nbw}

\numberwithin{condition}{\nbw}

\numberwithin{example}{\nbw}

\numberwithin{claim}{\nbw}

\numberwithin{remark}{\nbw}

\numberwithin{question}{\nbw}

\numberwithin{goal}{\nbw}

\numberwithin{fact}{\nbw}

\newcommand{\thmref}[1]{Theorem~\ref{thm:#1}} 
\newcommand{\propref}[1]{Proposition~\ref{prop:#1}} 
\newcommand{\corref}[1]{Corollary~\ref{cor:#1}} 
\newcommand{\secref}[1]{Section~\ref{sec:#1}} 
\newcommand{\eqnref}[1]{(\ref{eq:#1})} 

\def\be{\begin{equation} }
\def\ee{ \end{equation}}

\def\ben{\begin{equation*}}
\def\een{\end{equation*}}
\def\bea{\begin{eqnarray}}
\def\eea{\end{eqnarray}}
\def\ee{\end{eqnarray}}
\def\bean{\begin{eqnarray*}}
\def\eean{\end{eqnarray*}}

\newcommand\ignore[1]{}

\def\R{\mathbb{R}} 
\def\C{\mathbb{C}} 
\def\N{\mathbb{N}} 

\newcommand{\Ex}[1]{\mathbb{E}\left[#1\right]} 
\newcommand{\Prwo}{\mathbb{P}} 

\newcommand{\Var}[1]{\mathbb{V}\left(#1\right)} 
\renewcommand{\Pr}[1]{\mathbb{P}\left(#1\right)} 

\newcommand{\bigoh}[1]{O\left(#1\right)}
\newcommand{\liloh}[1]{o\left(#1\right)}
\newcommand{\ohmega}[1]{\Omega\left(#1\right)}

\def\sF{\mathcal{F}}

\def\sL{\mathcal{L}}
\def\sM{\mathcal{M}}

\def\sS{\mathcal{S}}

\def\deg{{\rm deg}}

\newcommand\QED{\ifhmode\allowbreak\else\nobreak\fi
\quad\nobreak$\Box$\medbreak}
\newcommand{\proofstart}{\par\noindent\sl Proof:\rm\enspace}
\newcommand{\proofend}{\QED\par}
\newenvironment{proof}{\proofstart}{\proofend}

\def\eps{\epsilon}

\def\deg{\indeg}

\renewcommand{\deg}{\mbox{d}}

\def\Id{I}

\def\cb{{\bf e}}

\begin{document}

\title{The spectrum of random $k$-lifts of large graphs (with possibly large $k$)}
\author{Roberto Imbuzeiro Oliveira\thanks{IMPA, Rio de Janeiro, RJ,
Brazil, 22430-040. \texttt{rimfo@impa.br}}} \maketitle

\begin{abstract}\ignore{A $k$-lift of a graph $G$ is a new graph
$G^{(k)}$ where each vertex $v$ of $G^{(k)}$ is replaced by $k$
copies itselft and each edge $vw$ is replaced by a matching between
the copies of $v$ and $w$.} We study random $k$-lifts of large, but
otherwise arbitrary graphs $G$. We prove that, with high
probability, all eigenvalues of the adjacency matrix of the lift
that are not eigenvalues of $G$ are of the order
$\bigoh{\sqrt{\Delta \ln (kn)}}$, where $\Delta$ is the maximum
degree of $G$. Similarly, and also with high probability, the ``new"
eigenvalues of the Laplacian of the lift are all in an interval of
length $\bigoh{\sqrt{\ln (nk)/d}}$ around $1$, where $d$ is the
minimum degree of $G$.

We also prove that, from the point of view of Spectral Graph Theory,
there is very little difference between a random $k_1k_2\dots
k_r$-lift of a graph and a random $k_1$-lift of a random $k_2$-lift
of $\dots $ of a random $k_r$-lift of the same graph.

The main proof tool is a concentration inequality for sums of random
matrices that was recently introduced by the author.\end{abstract}

\section{Introduction}

Let $G$ be a graph with vertex set $V$ and edge set $E$. A {\em
$k$-lift} of $G$ is a graph $G^{(k)}$ with vertex set $V\times [k]$
and edge set:$$E^{(k)}\equiv \cup_{vw\in E} \sM_{vw}$$ where each
$\sM_{vw}$ is a matching of the sets $\{(v,1),(v,2),\dots,(v,k)\}$ and
$\{(w,1),(w,2),\dots,(w,k)\}$. In more intuitive terms: each vertex
of $G$ is replaced by $k$ copies of itself and each edge $vw\in E$
is replaced by a matching of the copies of $v$ and $w$.

There have been many recent results about {\em random} $k$-lifts of
graphs where $G$ is fixed and $k\to +\infty$. Here ``random" means
that the matchings $\sM_{vw}$ are chosen independently and each of
them is uniformly distributed. A lot is now known about properties
of $G^{(k)}$ such as connectivity
\cite{AmitLinial_RandomLiftsIntro,AmitLinial_RandomLiftsEdgeExpansion},
chromatic number \cite{AmitLinialMatousek_RandomLiftsChromatic},
spectral distribution
\cite{Friedman_RelativelyRamanujan,LinialPuder_SpectraLifts} and the
existence of perfect matchings
\cite{LinialRozenman_MatchingsInLifts}.

A disjoint line of work has considered $2$-lifts of arbitrary
(possibly large) graphs $G$. The goal in this case was to provide an
explicit construction of {\em some} $2$-lift with good spectral
properties, so that arbitrarily large expanders can be efficiently
constructed via successive $2$-lifts \cite{BiluLinial_TwoLifts}.

In this paper we study a scenario that is quite natural but, to the
best of our knowledge, new: random $k$-lifts of large graphs $G$. We
obtain non-trivial results only when the minimum degree of $G$ is
$\gg \ln(|V|k)$, but $G$ and $k$ are otherwise arbitrary. For
concrete examples, one may think of random $n$-lifts of graphs on
$n$ vertices and minimal degree $\ln^{1+\eps}n$; or of
$2^{\sqrt{n}}$-lifts of $(n/2)$-regular graphs on $n$ vertices.

Our focus will be on the spectra of the {\em adjacency matrix} and
{\em Laplacian} of the random lift. These two matrices are the
central objects of Spectral Graph Theory and their eigenvalues can
be used to estimate many parameters of graphs, including the
diameter, distances between distinct subsets, discrepancy-like
properties, path congestion, cuts, chromatic number and the mixing
time for random walk; see e.g.
\cite{Chung_SpectralGraphTheory,ChungGraham_QuasirandomGraphsGivenDegrees,ChungGrahamWilson_QuasirandomGraphs}.
Our main theorem is a first indication of what the above parameters
are for the random lifts we consider. In fact, our theorem works
even for a relaxed definition of random lifts where the $\sM_{vw}$
need not be uniformly distributed.

We first need some preliminaries. Let $A$ and $A^{(k)}$ be the
adjacency matrix of the graph $G$ and of its $k$-lift $G^{(k)}$
(resp.). We will see in \secref{tensorlift} that the spectrum of $A^{(k)}$
always contains the spectrum of $A$ in the sense of {\em multisets}:
any eigenvalue of $A$ with multiplicity $m$ is an eigenvalue of
$A^{(k)}$ with multiplicity $\geq m$. The same holds for the spectra
of the Laplacians $\sL^{(k)}$ and $\sL$ of $G^{(k)}$ and $G$
(respectively).

Let ${\rm new}(A^{(k)})$ be the difference between the spectrum of
$A^{(k)}$ and the spectrum of $A$ and define ${\rm new}(\sL^{(k)})$
similarly. ${\rm new}(A^{(k)})$ is also a multiset: if $\lambda$ has
multiplicity $m_1$ in the spectrum of $A$ and multiplicity $m_2$ in
the spectrum of $A^{(k)}$, it occurs $m_2-m_1$ times in ${\rm
new}(A^{(k)})$. Our main result is:

\begin{theorem}\label{thm:randomlifts}With the above notation, let $n=|V|$ be the number of vertices in
$G$. Also let $d$ and $\Delta$ be the minimum and maximum degrees in
$G$ (respectively). Assume that the matchings $\{\sM_{vw}\}_{vw\in E}$ are chosen
independently and that for each $vw\in E$ and $\ell,r\in [k]$:
$$\Pr{\{(v,\ell),(w,r)\}\in\sM_{vw}}=\frac{1}{k}.$$
Then for all $\delta\in(0,1)$,
$$\Pr{\sup_{\eta\in {\rm new}(A^{(k)})}|\eta|\leq
16\sqrt{\Delta\ln(2nk/\delta)}}\geq 1-\delta$$ and
$$\Pr{\sup_{\beta\in {\rm new}(\sL^{(k)})}|1-\beta|\leq
16\sqrt{\frac{\ln(2nk/\delta)}{d}}}\geq 1-\delta.$$\end{theorem}

This is interesting even in the case $k=2$. It is known
\cite{BiluLinial_TwoLifts} that any $d$-regular graph has a two-lift
whose new eigenvalues are all $\bigoh{\sqrt{d\ln d}}$. However, a
{\em typical} random $2$-lift of $G_n$ might have at least one
eigenvalue equal to $d$. One example (also from
\cite{BiluLinial_TwoLifts}) consists of $n/(d+1)$ disconnected
$(d+1)$-cliques; the new eigenvalue $d$ comes from there being a
clique whose lift consists of two disconnected cliques.  [It is
possible to find connected examples with similar behavior.] Notice
that the probability of there being such a clique is $1-\liloh{1}$
even when $d=\lceil c\sqrt{\ln n}\rceil$ for some small constant
$c>0$. On the other hand, the Theorem shows that there exists some
$C>0$ such that for any $\eps>0$, if $d\geq C\ln n/\eps^2$, then the
largest new eigenvalue is $\leq \eps d$ with probability $\geq
1-1/n^2$.

On the other hand, we note that the largest eigenvalue of $A^{(k)}$
is always between $d$ and $\Delta$ and the eigenvalues of
$\sL^{(k)}$ are always between $0$ and $2$
\cite{Chung_SpectralGraphTheory}. Hence our result for the adjacency
matrix is trivial if $\Delta\leq \ln(nk/\delta)$ and the bound for
the Laplacian is trivial when $d\leq \ln(nk/\delta)$.

One corollary of \thmref{randomlifts} is the following result.

\begin{corollary}\label{cor:veryclose}In the setting of \thmref{randomlifts}, let $k=k_1\dots k_s$ with $k_1,\dots,k_s\in\N\backslash\{0,1\}$ and consider two different random
graphs:
\begin{itemize}
\item ${G}^{(k)}$ is a {\em maximally random} $k$-lift of $G$:
that is to say, each random matching $\sM_{vw}$ appearing in the
construction of $G^{(k)}$ is uniformly distributed over all
matchings of $\{(v,i)\}_{i=1}^k$ and $\{(w,j)\}_{j=1}^k$, and the
matchings are independent.
\item $\tilde{G}^{(k)}=G_s$ where $G_0=G$ and, for each $1\leq i\leq
s$, $G_i$ is a maximally random $k_i$-lift of $G_{i-1}$ (conditionally on $G_0,G_1,\dots,G_{i-1}$).
\end{itemize}
Let $A^{(k)}$ and $\sL^{(k)}$ denote the adjancency matrix and
Laplacian of $G^{(k)}$ and define $\tilde{A}^{(k)}$ and
$\tilde{\sL}^{(k)}$ similarly. Then (with an appropriate labelling
of the vertices of the two graphs):
$$\Pr{\|A^{(k)}-\tilde{A}^{(k)}\|\leq
32\,\sqrt{\Delta\ln(4nk/\delta)}}\geq 1-\delta$$ and
$$\Pr{\|\sL^{(k)}-\tilde{\sL}^{(k)}\|\leq
32\,\sqrt{\frac{\ln(4nk/\delta)}{d}}}\geq 1-\delta.$$\end{corollary}

This is interesting because the distributions of $G^{(k)}$ and $\tilde{G}^{(k)}$
can be very different. For instance, let $k_1=k_2=\dots=k_s=2$. If $s$ is constant and the number of vertices is large enough, all $2^s!$ possible permutations will be seen in the matchings of $\{v\}\times [k]$ with $\{w\}\times [k]$ for $vw\in E$. On the other hand, only $2^s$ possible permutations will be seen in $\tilde{G}^{(k)}$.

\thmref{randomlifts} will be deduced from a recent concentration
result for sums of independent random {\em matrices}. In what follows
$\C^{d\times d}_{\rm Herm}$ is the space of $d\times d$ Hermitian
matrices with complex entries, the expectations of matrices are
defined entrywise and $\|\cdot\|$ is the operator norm. [See
\secref{prelimlinearalgebra} and \secref{prelimprob} for these and
related definitions.]

\begin{theorem}[Corollary 7.1 in \cite{Eu_Freedman}]\label{thm:freedman}Let $X_1,\dots,X_m$ be mean-zero independent random matrices, defined on a common probability space $(\Omega,\sF,\Prwo)$, with values in $\C^{d\times d}_{\rm Herm}$ and such that there exists a $M>0$ with $\|X_i\|\leq M$ almost surely for all $1\leq i\leq m$. Define:
$$\sigma^2\equiv \mbox{ the largest eigenvalue of }\sum_{i=1}^m\Ex{X_i^2}.$$
Then for all $t\geq 0$,
$$\Pr{\left\|\sum_{i=1}^m X_i\right\|\geq t}\leq 2d\, e^{-\frac{t^2}{8\sigma^2 + 4Mt}}.$$\end{theorem}

Given this bound, \thmref{randomlifts} follows quite easily, while
other proof techniques for bounding spectra of random matrices (such
as the {\em trace method}
\cite{FurediKomlos_Eigenvalues,Friedman_RelativelyRamanujan,Friedman_SecondEigenvalue,LinialPuder_SpectraLifts}
and the discrepancy-based ideas of Feige and Ofek
\cite{FeigeOfek_Spectral}) can be quite technical. In our setting,
\thmref{freedman} is also an improvement over other general
concentration bounds for random matrices, most notably the operator
Chernoff bound of Ahlswede and Winter
\cite{AhlswedeWinter_StrongConverse} and the matrix Hoeffding bound
of Christofides and Markstr\"{o}m
\cite{ChristofidesMarkstrom_HoeffdingForMatrices}. A key advantage
of \thmref{freedman} over related results is that its ``variance"
term can be much smaller, especially in the graph-theoretical
setting; this is discussed in more detail in Remark 7.1 of
\cite{Eu_Freedman}.

The remainder of the paper is organized as follows. After the
preliminary \secref{prelim}, we collect some basic facts about
$k$-lifts in \secref{basiclifts}. We prove the Theorem and its
Corollary in \secref{mainproofs}. The last Section presents some
extensions and open questions.\\

\section{Preliminaries}\label{sec:prelim}

\subsection{Basic notation}

For a natural number $m\in\N\backslash\{0\}$, $[m]$ is the set of
all integers $1\leq i\leq m$.

We will frequently speak of {\em multisets} $S$. Given a ground set
$\sS$ (which will usually be $\R$), a multiset $S$ is defined by a
function $m_S:\sS\to\N$. Informally, we will let think of $S$ as a
set where each $x\in A$ appears $m_S(x)$ times and we will refer to
this quantity as the {\em multiplicity} of $x$. We say that $x$
belongs to $S$ ($x\in S$) if $m_S(x)>0$.

For two multisets $S_1,S_2$ over the same ground set $\sS$ and with
corresponding functions $m_{S_1},m_{S_2}$, we say that $S_1\subset
S_2$ if for all $x\in \sS$ $m_{S_1}(x)\leq m_{S_2}(x)$. The
difference $S_2\backslash S_1$ is the multiset where each $x\in \sS$
has multiplicity $\max\{m_{S_2}(x)-m_{S_1}(x),0\}$.

\subsection{Linear algebra}\label{sec:prelimlinearalgebra}

For given $d_r,d_c\in\N\backslash\{0\}$, $\R^{d_r\times d_c}$ (resp.
$\C^{d_r\times d_c}$) is the space of $d_r\times d_c$ matrices with
entries in $\R$ (resp. $\C$).

For $A\in\R^{d_r\times d_c}$, $A^\dag\in\R^{d_c\times d_r}$ is the
transpose of $A$; similarly, for $B\in \C^{d_r\times d_c}$,
$B^*\in\C^{d_c\times d_r}$ is the conjugate transpose of $B$. We
identify $\R^d$ and $\C^d$ with $\R^{d\times 1}$ and $\C^{d\times
1}$ (resp.), so that the standard inner product of $x,y\in \R^d$ is
$x^\dag y$.

$\C^{d\times d}_{\rm Herm}$ is the space of $d\times d$ {\em
Hermitian matrices}, which are the $A\in\C^{d\times d}$ with
$A^*=A$. Similarly, $\R^{d\times d}_{\rm Sym}$ is the space of all
$d\times d$ real matrices that are symmetric in the sense that
$A=A^\dag$.

For a vector $v\in\R^d$ or $\C^d$, $\|v\|$ is its Euclidean norm. The {\em operator norm} of $A\in \R^{d\times d}$ is:
$$\|A\|\equiv \max_{v\in\R^d,\, \|v\|=1}\|Av\|.$$

Finally, the canonical basis vectors for $\R^d$ is
denoted by $\cb_1,\cb_2,\dots,\cb_d$.

\subsubsection{The spectral theorem}

We recall the standard spectral theorem: for any $A\in \R^{d\times
d}_{\rm Sym}$ there exists a set $S\subset \R$ and orthogonal
projections $\{P_\alpha\}_{\alpha\in S}$ with orthogonal ranges such
that:
$$\sum_{\alpha\in S}\alpha \,P_\alpha=A\mbox{ and }\sum_{\alpha\in
S}P_\alpha = I_d,$$ where $I_d$ is the $d\times d$ identity matrix.
The numbers $\alpha\in S$ are called the eigenvalues of $A$ and the
vectors $v$ in the range of $P_\alpha$ are eigenvectors
corresponding to a given $\alpha$. The spectrum of $A$, denoted by
${\rm spec}(A)$, is the multiset where each $\alpha\in S$ appears
with multiplicity equal to the rank of $P_\alpha$.

One useful consequence of the spectral decomposition is that $\|A\| = \max_{\alpha\in{\rm spec}(A)}|\alpha|$.

\subsubsection{Tensor products}\label{sec:prelimtensor}

It will be convenient to represent the matrices of lifts via tensor
products. The {\em tensor product} of $\R^{d_1}$ and $\R^{d_2}$,
denoted by $\R^{d_1}\otimes \R^{d_2}$, is the set of all formal
linear combinations of vectors of the form $\cb_{i_1}\otimes
\cb_{i_2}$ with $1\leq i_b\leq d_b$ for $b=1,2$. [We will abuse
notation and assume that $e_i\in \R^{d_1}\cap \R^{d_2}$ for $i\leq
\min\{d_1,d_2\}$.]

Similarly, if $v_{b}=\sum_{j_b=1}^{d_b}v_{b,j_b}\cb_{j_b}$
($b=1,2$), the tensor product of $v_1\otimes v_2$ is defined by the
``distributive rule":
$$v_1\otimes v_2\equiv
\sum_{j_1=1}^{d_1}\sum_{j_2=1}^{d_2}v_{1,j_1}v_{2,j_2}\cb_{j_1}\otimes
\cb_{j_2}.$$

There exists a unique inner product on $\R^{d_1}\otimes \R^{d_2}$,
denoted by $(\cdot,\cdot\cdot)$, such that for all
$v_1,w_1\in\R^{d_1}$ and $v_2,w_2\in\R^{d_2}$,
$$(v_1\otimes v_2,w_1\otimes w_2)=(v_1^\dag w_1)\,(v_2^\dag w_2).$$

Moreover, the tensor product of $A_{1}\in\R^{d_1\times d_1}$ and
$A_2\in\R^{d_2\times d_2}$ is the unique linear operator $A_1\otimes
A_2$ from $\R^{d_1}\otimes \R^{d_2}$ to itself that satisfies:
$$\forall 1\leq i_1\leq d_1,\,\forall 1\leq i_2\leq d_2,\; (A_1\otimes
A_2)(\cb_{i_1}\otimes \cb_{i_2}) = (A_1\cb_1)\otimes (A_2\cb_2).$$

One can check that if $A_1\in \R_{\rm Sym}^{d_1\times d_1}$ and
$A_2\in \R_{\rm Sym}^{d_2\times d_2}$, then $A_1\otimes A_2$ is
self-adjoint in the sense that:

$$\forall u,v\in \R^{d_1}\otimes \R^{d_2},\, (u,(A_1\otimes A_2)v) = ((A_1\otimes A_2)u,v).$$

In general, one still has:

\begin{equation}\label{eq:tensoradjoint}\forall u,v\in \R^{d_1}\otimes \R^{d_2},\, (u,(A_1\otimes
A_2)v) = ((A_1^\dag\otimes A_2^\dag)u,v).\end{equation}

i.e. $A_1^\dag \otimes A_2^\dag$ is the adjoint of $A_1\otimes A_2$.

Notice that $\R^{d_1}\otimes \R^{d_2}$ is isomorphic to
$\R^{d_1d_2}$, in the sense that any bijection $\psi:[d_1]\times
[d_2]\to [d_1d_2]$ can be ``lifted" to an invertible,
inner-product-preserving linear map:
$$\Psi: \R^{d_1}\otimes \R^{d_2} \to \R^{d_1d_2}$$
defined by the rule $\Psi(\cb_i\otimes\cb_j) = \cb_{\psi(i,j)}$,
$(i,j)\in[d_1]\times [d_2]$. Under this map, self-adjoint maps over
$\R^{d_1}\otimes\R^{d_2}$ correspond to symmetric matrices over
$\R^{d_1d_2}$ and vice versa. Therefore, one may also state a
spectral theorem over $\R^{d_1}\otimes \R^{d_2}$; we omit the
details.

\subsection{Concepts from Graph
Theory}\label{sec:prelimgraph}

For our purposes a graph $G=(V,E)$ consists of a finite set $V$ of
vertices and a set $E$ of edges, which are subsets of size $2$ of
$V$. Unless otherwise noted, we will assume that $V=[n]$ for some
integer $n\geq 2$, where $[n]\equiv \{1,2,\dots,n\}$. We will write
edges as unordered pairs $vw$ or $\{v,w\}$ and make no distinction
between $vw$ and $wv$. The {\em degree} $\deg_G(v)$ of a vertex $v$
is the number of $w\in V\backslash\{v\}$ such that $vw\in E$.

Assume that $V=[n]$, or more generally, that the elements of $V$ are
labelled $v_1,\dots,v_n$. The {\em adjacency matrix} of $G$ is the
$n\times n$ matrix $A\in\R^{n\times n}_{\rm Sym}$ with zeros on
the diagonal and such that, for all $1\leq i<j\leq n$, the
$(i,j)$-th entry of $A$ is $1$ if $v_iv_j\in E$ and $0$ otherwise.
When $V=[n]$, this reads: \begin{equation}\label{eq:whatAis}A\equiv
\sum_{ij\in E}(\cb_i\cb_j^\dag +\cb_j\cb_i^\dag).\end{equation}

The {\em Laplacian} $\sL$ of $G$ is the matrix:
$$\sL = \Id_n - T\,A\,T$$
where $T$ is the $n\times n$ diagonal matrix whose $(i,i)$-th entry
is $\deg_G(i)^{-1/2}$ if $\deg_G(i)\neq 0$, or $0$ if $\deg_G(i)=0$.
If all degrees are non-zero, one can write this as follows:
\begin{equation}\label{eq:whatsLis}\sL = \Id_n - \sum_{ij\in E}\frac{(\cb_i\cb_j^\dag
+\cb_j\cb_i^\dag)}{\sqrt{\deg_G(i)\deg_G(j)}}.\end{equation}
\ignore{ We also let
$$\lambda(G)\equiv
\min\{\lambda_1(\sL),2-\lambda_{d-1}(\sL)\}$$ denote the {\em
spectral gap} of $G$.}

\subsection{Probability with
matrices}\label{sec:prelimprob}

We will be dealing with random matrices (and random linear
operators) throughout the paper. Following common practice, we will
always assume that we have a probability space $(\Omega,\sF,\Prwo)$
in the background where all random variables are defined.

Call a map $X:\Omega\to\C^{d\times d}_{\rm Herm}$ a random $d\times
d$ Hermitian matrix (or a $\C^{d\times d}_{\rm Herm}$-valued random
variable) if for each $1\leq i,j\leq n$, the function
$X(i,j):\Omega\to \C^{d\times d}$ corresponding to the $(i,j)$-th
entry of $X$ is $\sF$-measurable, or equivalently, if for each Borel
subset $S\subset \C^{d\times d}_{\rm Herm}$ $X^{-1}(S)\in\sF$. We
say that $X$ is integrable if all the entries of $X$ are integrable,
one defines $\Ex{X}$ entrywise: the $(i,j)$th entry of $\Ex{X}$ is
$\Ex{X(i,j)}$. We will also use analogous definitions for
$X:\Omega\to\R^{d\times d}_{\rm Sym}$. [We will essentially ignore
all measurability and integrability issues in the remainder of the
paper. These can be dealt with in a rather straightforward manner.]

One can easily check that if $X$ is a random integrable $d\times d$
Hermitian matrix and $A\in\C^{d\times d}_{\rm Herm}$ is
deterministic, $\Ex{AX} = A\Ex{X}.$ If the entries
of $X$ are also square integrable, one may define a ``matrix
variance" $\Var{X} \equiv \Ex{(X-\Ex{X})^2}$ and deduce that:
\begin{equation}\label{eq:matrixvariance}\Var{X} = \Ex{X^2} -
\Ex{X}^2.\end{equation}

\section{Lifts of graphs}\label{sec:basiclifts}

Our goal here is to review the construction of lifts of graphs
outlined in the introduction and to prove some elementary facts that
will be useful later on. Other perspectives on these objects can be
found in \cite{AmitLinial_RandomLiftsIntro}.

Recall that a {\em matching} of two finite, disjoint, non-empty sets
$A,B$ is a set of pairs:
$$\sM=\{\{a_i,b_i\}\,:\,i=1,\dots,m\}$$
where $(a_1,\dots,a_m)$ is a permutation of the elements of $A$ and
$(b_1,\dots,b_m)$ is a permutation of the elements of $B$. Notice
that the existence of a matching $\sM$ as above implies that
$|A|=|B|=m$.

Now let $G$ be a graph with vertex set $V=[n]$ and edge set $E$.
Given $k\in\N\backslash\{0,1\}$, a $k$-lift $G^{(k)}$ of $G$ is
determined by a choice of matchings: $$\{\sM_{ij}\,:\,ij\in E \},$$
where for each $ij\in E$ $\sM_{ij}$ is a matching of $\{i\}\times
[k]$ and $\{j\}\times [k]$. $G^{(k)}$ is the graph with vertex set
$[n]\times [k]$ and edge set $E^{(k)}=\cup_{ij\in E}\sM_{ij}$.

\subsection{Graph matrices and tensor
products}\label{sec:tensorlift}

It is convenient to represent the matrices corresponding to
$G^{(k)}$ in the tensor space $\R^n\otimes \R^k$. That is to say, we
will write down a linear operator $A^{(k)}$ over $\R^n\otimes \R^k$
such that for all $(i,\ell),(j,r)\in [n]\times [k]$, $$(\cb_i\otimes
\cb_\ell, A^{(k)}(\cb_j\otimes\cb_r)) = \left\{\begin{array}{ll}1 &
\mbox{if } \{(i,\ell),(j,r)\}\in E^{(k)};\\ 0 &
\mbox{otherwise}.\end{array}\right.$$ This is satisfied by:
$$A^{(k)} = \sum_{\{(i,\ell),(j,r)\}\in E^{(k)}}(\cb_{i}\cb_j^\dag
)\otimes (\cb_\ell\cb_r^\dag) + (\cb_{j}\cb_i^\dag )\otimes
(\cb_r\cb_\ell^\dag).$$

Another way of writing $A^{(k)}$ will be more useful later on:
\begin{eqnarray}\label{eq:rewriteAk}A^{(k)} &=& \sum_{ij\in
E}\cb_{i}\cb_j^\dag\otimes V_{(i,j)} + \cb_{j}\cb_i^\dag
\otimes V_{(j,i)},\mbox{ where $V_{(i,j)}$ is defined as:}\\
\label{eq:defMij} V_{(i,j)} &\equiv& \sum\limits_{(\ell,r)\in
[k]^2\,:\, \{(i,\ell),(j,r)\}\in \sM_{ij}}
\cb_\ell\cb_r^\dag.\end{eqnarray} We emphasize that the definition
of $V_{(i,j)}$ is {\em not} symmetric with respect to $i,j$: in
fact, a simple computation shows that $V_{(j,i)}=V_{(i,j)}^\dag =
V_{(i,j)}^{-1}$.

The Laplacian $\sL^{(k)}$ of $G^{(k)}$ can be similarly written as a
linear operator over $\R^n\otimes \R^k$. The key point to notice is
that all copies of $i\in[n]$ in $G^{(k)}$ have the same degree,
i.e.:
$$\forall \ell\in[k],\, \deg_{G^{(k)}}((i,\ell)) = \deg_{G}(i).$$
A simple calculation (omitted) shows that:
\begin{equation}\label{eq:rewriteLk}\sL^{(k)} = I_n\otimes I_k - \sum_{ij\in E}\frac{\cb_{i}\cb_j^\dag\otimes V_{(i,j)} + \cb_{j}\cb_i^\dag
\otimes V_{(j,i)}}{\sqrt{\deg_{G}(i)\deg_{G}(j)}}.\end{equation}

\subsection{Old and new eigenvalues}

We now draw a connection between the spectrum and eigenvalues of $A$ and
$A^{(k)}$. All arguments here also appear on previous papers on
graph lifts (e.g.
\cite{BiluLinial_TwoLifts,Friedman_RelativelyRamanujan}).

\begin{proposition}\label{prop:new}The spectrum of the adjacency matrix $A$ of $G$ is contained in the spectrum
of $A^{(k)}$ (counting multiplicities). Moreover, if $${\rm
new}(A^{(k)})\equiv{\rm spec}(A^{(k)})\backslash {\rm spec}(A)$$ is
the difference of the two spectra as multisets,
$$\max_{\eta\in {\rm new}(A^{(k)})}|\eta|=\|A^{(k)} - A\otimes \Pi_k\|$$ where $\Pi_k$ is the
$k\times k$ matrix with all entries equal to $1/k$.\end{proposition}

Essentially the same argument shows a related result for the
Laplacian $\sL^{(k)}$ of $G^{(k)}$ (proof omitted).

\begin{proposition}\label{prop:newsL}The spectrum of the Laplacian $\sL$ of $G$ is contained in the spectrum
of $\sL^{(k)}$ (counting multiplicities). Moreover, if $${\rm
new}(\sL^{(k)})\equiv{\rm spec}(\sL^{(k)})\backslash {\rm
spec}(\sL)$$ is the difference of the two spectra as multisets,
$$\max_{\eta\in {\rm new}(\sL^{(k)})}|1-\eta|=\|\sL^{(k)} - (I_n\otimes I_k - (I-\sL)\otimes \Pi_k)\|$$ where $\Pi_k$ is the
$k\times k$ matrix with all entries equal to $1/k$.\end{proposition}

\begin{proof}[of \propref{new}] Let ${\bf 1}_k\in\R^k$ be the vector with all coordinates equal to
$1$. Notice that $V_{(i,j)}{\bf 1}_k = \Pi_k{\bf 1}_k = {\bf 1}_k$
for all $i,j$ with $ij\in E$. Therefore, for all vectors $v\in
\R^n$,
$$A^{(k)}(v\otimes {\bf 1}_k) = (Av)\otimes {\bf
1}_k=(A\otimes \Pi_k) {\bf 1}_k.$$

In particular, if $v$ is an eigenvector of $A$ with eigenvalue
$\lambda$, $v\otimes {\bf 1}_k$ is an eigenvector of both $A^{(k)}$
and $A\otimes \Pi_k$, with the same eigenvalue $\lambda$ for both
matrices. It follows that each eigenvalue $\lambda$ of $A$ with
multiplicity $m$ is an eigenvalue of $A^{(k)}$ with multiplicity
$\geq m$, which is the first assertion in the Proposition.

Any new eigenvalue $\eta\in{\rm new}(A^{(k)})$ must correspond an
eigenvector $w\in\R^n\otimes \R^k$ that is orthogonal to $v\otimes
{\bf 1}_k$ for all eigenvectors $v$ of $A$ corresponding to ``old"
eigenvalues. Since the eigenvectors of $A$ span $\R^n$, any $w$ as
above must be orthogonal to the subspace:
$$H\equiv \{v\otimes {\bf 1}_k\,:\, v\in\R^n\}\subset\R^n\otimes\R^k.$$

In particular, $$\max_{\eta\in{\rm new}(A^{(k)})}|\eta| = \max_{w\in
H^\perp}\|A^{(k)}w\|.$$

To finish, we must show that the RHS equals $\|A^{(k)}-A\otimes
\Pi_k\|$. We have already seen that the operators $A^{(k)}$ and
$A\otimes \Pi_k$ have $H$ as an invariant subspace and that their
restrictions to that subspace are equal. This implies that $H^\perp$
must also be invariant and moreover:
$$\|A^{(k)} - A\otimes \Pi_k\| = \max_{w\in H^\perp\,:\, \|w\|=1} \|A^{(k)}w - (A\otimes
\Pi_k)w\|.$$ Now notice that:
$$H^\perp\equiv {\rm span}\{x\otimes y\,:\, x\in\R^n, y\in\R^k,\,y\perp {\bf
1}_k\}.$$ Moreover, for all $x\otimes y$ as above,
$$(A\otimes \Pi_k)(x\otimes y) = (Ax)\otimes (\Pi_ky) = 0$$
since $\Pi_k$ is the projection onto the line spanned by ${\bf
1}_k$. By linearity, this implies that $(A\otimes \Pi_k) w=0$ for
all $w\in H^\perp$, which results in the desired equality:
$$\|A^{(k)} - A\otimes \Pi_k\| = \max_{w\in H^\perp\,:\, \|w\|=1} \|A^{(k)}w\|.$$
\end{proof}

\section{Main proofs}\label{sec:mainproofs}

Propositions \ref{prop:new} and \ref{prop:newsL} show that in order
to prove \thmref{randomlifts}, one must bound the difference between
certain matrices. We attack this problem from the perspective of
concentration of measure. As it turns out, $A\otimes \Pi_k$ is the
expected value of $A^{(k)}$ and $I_n\otimes I_k - (I_n-\sL)\otimes
\Pi_k$ is the expected value of $\sL^{(k)}$. The concentration
inequality in \thmref{freedman} will ensure that $A^{(k)}$ and
$\sL^{(k)}$ are likely to be close to their respective expected
values. One this is achieved, \thmref{randomlifts} and its Corollary
will easily follow.

\subsection{Proof of the main theorem}

In this section we prove \thmref{randomlifts}.

\begin{proof}[of \thmref{randomlifts}] We start with the result for the adjacency matrix. \propref{new} implies that it is necessary and sufficient to prove that:
\begin{equation}\label{eq:tudoestaaqui}{\bf [Goal]}\;\;\Pr{\|A^{(k)}-A\otimes \Pi_k\|\leq 16\sqrt{\Delta\ln(2nk/\delta)}}\geq 1-\delta.\end{equation}
We will restate this as a concentration bound for the sum of random
matrices. Recall from \secref{tensorlift} that:
$$A^{(k)} = \sum_{ij\in E}Z_{ij} \mbox{ where }Z_{ij} =
\cb_j\cb_i^\dag \otimes V_{(i,j)} + \cb_i\cb_j^\dag \otimes
V_{(j,i)}.$$ We notice that all $Z_{ij}$ are self-adjoint, as
attested by \eqnref{tensoradjoint} and the fact that $V_{(i,j)}^\dag
= V_{(j,i)}$ (cf. \secref{tensorlift}).

The matrices $V_{(i,j)}$ and $V_{(j,i)}$ are determined by the
random matching $\sM_{ij}$. Since these matchings are independent,
the $\{Z_{ij}\}_{ij\in E}$ are also independent. Let us now compute
$\Ex{Z_{ij}}$ for a fixed $ij\in E$. It is not hard to show that
this is:
$$\Ex{Z_{ij}}=\cb_j\cb_i^\dag\otimes \Ex{V_{(i,j)}} + \cb_i\cb_j^\dag \otimes
\Ex{V_{(j,i)}}.$$ The $(\ell,r)$-th entry of $V_{(i,j)}$ is an
indicator random variable that is equal to $1$ iff $(i,\ell)$ and
$(j,r)$ are connected in the matching. By assumption, this happens
with probability $1/k$, therefore each entry of $V_{(i,j)}$ has
expected value $1/k$. This implies that $\Ex{V_{(i,j)}}$ is
precisely the matrix $\Pi_k$ in the Theorem. Similarly,
$\Ex{V_{(j,i)}}=\Pi_k$. We deduce that:
\begin{equation}\label{eq:expectedZij}\Ex{Z_{ij}} = (\cb_j\cb_i^\dag +
\cb_i\cb_j^\dag)\otimes \Pi_k.\end{equation} Now employ
\eqnref{whatAis} to deduce that:
$$\sum_{ij\in E}\Ex{Z_{ij}} = \left(\sum_{ij\in E}\cb_j\cb_i^\dag +
\cb_i\cb_j^\dag\right)\otimes \Pi_k. = A\otimes \Pi_k.$$ In other
words,
\begin{equation}\label{eq:tudoagora}A^{(k)}-A\otimes \Pi_k = \sum_{ij\in E}(Z_{ij}-\Ex{Z_{ij}})\end{equation}
is a sum of independent, self-adjoint random linear operators with
mean $0$. One may recall from \secref{prelimtensor} that
self-adjoint linear operators over $\R^{n}\otimes\R^k$ correspond to
symmetric matrices over $\R^{nk}$; therefore, we can apply
\thmref{freedman} to the above sum once we compute the variance
parameter $\sigma^2$ and the uniform bound $M$.

We start with $M$. $Z_{ij}$ is the adjacency matrix of a graph that
has all degrees equal to $1$. Therefore, $\|Z_{ij}\|=1$ and (by
Jensen's inequality) $\Ex{\|Z_{ij}\|}\leq 1$. It follows that all
terms in the sum \eqnref{tudoagora} satisfy
$\|Z_{ij}-\Ex{Z_{ij}}\|\leq M\equiv 2$.

To compute $\sigma^2$, we start with $\Ex{Z_{ij}^2}$ for a fixed
$ij\in E$. One can check that:
$$Z_{ij}^2 = \cb_i\cb_i^\dag\otimes (V_{(j,i)}V_{(i,j)}) + \cb_j\cb_j^\dag\otimes
(V_{(i,j)}V_{(j,i)}).$$ Now recall from \secref{tensorlift} that
$V_{(j,i)}=V_{(i,j)}^{-1}$ and deduce that:
$$Z_{ij}^2 = (\cb_i\cb_i^\dag + \cb_j\cb_j^\dag)\otimes I_k.$$
Another computation reveals that:
$$\Ex{Z_{ij}}^2 = (\cb_i\cb_i^\dag + \cb_j\cb_j^\dag)\otimes \Pi_k.$$
Using \eqnref{matrixvariance}, we deduce that:
$$\Var{Z_{ij}} = \Ex{(Z_{ij}-\Ex{Z_{ij}})^2} = (\cb_i\cb_i^\dag +
\cb_j\cb_j^\dag)\otimes (I_k - \Pi_k).$$ Summing up those terms, we
arrive at:
$$\sum_{ij\in E}\Ex{(Z_{ij}-\Ex{Z_{ij}})^2} = \left[\sum_{ij\in E} (\cb_i\cb_i^\dag +
\cb_j\cb_j^\dag)\right]\otimes (I_k - \Pi_k) = [\sum_{i=1}^n
\deg_G(i)\, \cb_i\cb_i^\dag]\otimes (I_k - \Pi_k).$$ Given two
symmetric matrices $B_1,B_2$, the eigenvalues of $B_1\otimes B_2$
are precisely the products of the form $\lambda_1\lambda_2$ with
$\lambda_i\in{\rm spec}(B_i)$, $i=1,2$. To apply this above, notice
that $\Pi_k$ is a rank-$1$ projection, hence the eigenvalues of
$I_k-\Pi_k$ are $0$ and $1$. It follows that:
$$\left\|[\sum_{i=1}^n
\deg_G(i)\, \cb_i\cb_i^\dag]\otimes (I_k - \Pi_k)\right\| =
\left\|\sum_{i=1}^n \deg_G(i)\, \cb_i\cb_i^\dag\right\|.$$ But the
matrix on the RHS is diagonal with non-negative entries, hence its
largest eigenvalue is the largest entry on the diagonal, which is
$\max_{i}\deg_G(i)=\Delta$. We deduce that one may take
$\sigma^2=\Delta$.

We now apply \thmref{freedman} with $\sigma^2=\Delta$ and $M=2$ to
the sum of the independent random linear operators
$\{Z_{ij}-\Ex{Z_{ij}}\}_{ij\in E}$, which is $A^{(k)}-A\otimes
\Pi_k$ (cf. \eqnref{tudoagora}). Moreover, the dimension parameter
in this case is $d=nk$ because that is the dimension of the space
$\R^n\otimes\R^k$ where the matrices are defined. We obtain:
$$\Pr{\|A^{(k)}-A\otimes \Pi_k\|\geq t}\leq 2nk
\,e^{-\frac{t^2}{8(\Delta+t)}}.$$ Taking $t \equiv
16\max\{\sqrt{\Delta \ln(2nk/\delta)},\ln(2nk/\delta)\}$ makes the
RHS smaller than $\delta$. This implies the desired result if
$\Delta\geq \ln(2nk/\delta)$. However, notice that $\|A^{(k)}\|\leq
\Delta$, as $G^{(k)}$ is a graph of maximal degree $\Delta$; and
similarly, $\|A\otimes \Pi_k\|\leq \Delta$. Therefore, we have
$\|A^{(k)}-A\otimes \Pi_k\|\leq 2\Delta$ {\em always} and this
implies that we still have the postulated bound if $\Delta\leq \ln
(2nk/\delta)$, as in that case $16\sqrt{\Delta\ln(2nk/\delta)}\geq
16\Delta.$ This proves \eqnref{tudoestaaqui}, which (as seen above)
is equivalent to the desired assertion via \propref{new}.

The proof for the Laplacian is quite similar and we will present it
in less detail. We use \propref{newsL} in order to restate the
desired inequality as:

\begin{equation}\label{eq:tudoestaaquisL}{\bf [Goal]}\;\;\Pr{\|\sL^{(k)}-(I_n\otimes I_k - (I-\sL)\otimes \Pi_k)\|\leq 16\sqrt{\frac{\ln(2nk/\delta)}{d}}}\geq 1-\delta.\end{equation}

Using equations \eqnref{whatsLis} and \eqnref{rewriteLk}, we see
that: \begin{eqnarray}\nonumber I_n\otimes I_k - (I-\sL)\otimes
\Pi_k - \sL^{(k)} &=& \sum_{ij\in E}\frac{\cb_j\cb_i^\dag \otimes
(V_{(i,j)}-\Pi_k) + \cb_i\cb_j^\dag \otimes (V_{(j,i)}
-\Pi_k)}{\sqrt{\deg_G(i)\deg_G(j)}}\\ \label{eq:tudoagorasL} &=&
\sum_{ij}\frac{Z_{ij}-\Ex{Z_{ij}}}{\sqrt{\deg_G(i)\deg_G(j)}}\end{eqnarray}
with the same $Z_{ij}$ from the first part. The terms in the sum are
again independent matrices with mean $0$ and we will apply
\thmref{freedman} to their sum. For this, we need to compute the
corresponding $M$ and $\sigma^2$.

For the parameter $M$, we observe that, since $d$ is the minimum
degree and $\|Z_{ij}-\Ex{Z_{ij}}\|\leq 2$ (as shown before),
$$\left\|\frac{Z_{ij}-\Ex{Z_{ij}}}{\sqrt{\deg_G(i)\deg_G(j)}}\right\|
\leq \frac{2}{\sqrt{\deg_G(i)\deg_G(j)}}\leq \frac{2}{d},$$ hence we
may take $M=2/d$. Each term in the sum has variance:
$$\Ex{\left(\frac{Z_{ij}-\Ex{Z_{ij}}}{\sqrt{\deg_G(i)\deg_G(j)}}\right)^2}
= \frac{1}{\deg_G(i)\deg_G(j)}\,\Ex{(Z_{ij}-\Ex{Z_{ij}})^2} =
\frac{(\cb_i\cb_i^\dag + \cb_j\cb_j^\dag)\otimes
(I_k-\Pi_k)}{\deg_G(i)\deg_G(j)}.$$ The sum of these terms is:
$$\sum_{ij\in E}\frac{(\cb_i\cb_i^\dag + \cb_j\cb_j^\dag)\otimes
(I_k-\Pi_k)}{\deg_G(i)\deg_G(j)} = \sum_{i=1}^n \left(\sum_{j:ij\in
E}\frac{\cb_i\cb_i^\dag}{\deg_G(i)\deg_G(j)}\right)\otimes
(I_k-\Pi_k).$$ Again we have a tensor product of a diagonal matrix
with another matrix whose eigenvalues are either $0$ or $1$. We
deduce as before that the operator norm is at most:
$$\max_{i}\sum_{j:ij\in
E}\frac{1}{\deg_G(i)\deg_G(j)}\leq \max_i \sum_{j:ij\in
E}\frac{1}{\deg_G(i)d} = \frac{1}{d}.$$ Therefore, we may take
$\sigma^2=1/d$.

Apply now \thmref{freedman} and \eqnref{tudoagorasL} to deduce that:
$$\Pr{\|\sL^{(k)}-(I_n\otimes I_k - (I-\sL)\otimes \Pi_k)\|\geq t}\leq
2nk\,e^{-\frac{t^2d}{8+8t}}.$$ Taking:
$$t\equiv
16\,\max\left\{\sqrt{\frac{\ln(2nk/\delta)}{d}},\frac{\ln(2nk/\delta)}{d}\right\}$$
makes the RHS $\leq \delta$ and implies the desired result when
$\ln(2nk/\delta)/d\leq 1$. However, any graph Laplacian has spectrum
contained in $[0,2]$ \cite{Chung_SpectralGraphTheory}; this implies
that $\|\sL^{(k)}-(I_n\otimes I_k - (I-\sL)\otimes \Pi_k)\|\leq 4$
always. In particular, the bound claimed in \eqnref{tudoestaaquisL}
holds even if $\ln(2nk/\delta)/d>1$. This finishes the proof of
\eqnref{tudoestaaquisL}, which implies the Theorem (cf.
\propref{newsL}).\end{proof}

\subsection{Proof of the corollary}

\begin{proof}[of \corref{veryclose}] We only present the proof of the adjacency matrix; the argument for the Laplacian is exactly the same.

The adjacency matrix $A^{(k)}$ of the graph $G^{(k)}$ satisfies:
\begin{equation}\label{eq:closepacas}\Pr{\|A^{(k)} - A\otimes \Pi^{(k)}\|\leq 16\sqrt{\Delta\ln(4nk/\delta)}}\geq 1 - \frac{\delta}{2}.\end{equation}
This is precisely what we showed in the course of the proof of \thmref{randomlifts} and also follows from applying the Theorem in conjunction with \propref{new}.

We claim that the same bound holds for $\tilde{A}^{(k)}$, after a suitable relabelling of the vertices. The vertex set of this graph  is $[n]\times K$ where $$K=[k_1]\times [k_2]\times \dots \times [k_s].$$ A simple induction argument shows that $\tilde{G}^{(k)}$ is also a lift of $G$, in the sense that its edge set $\tilde{E}^{(k)}$ is a union:
$$\tilde{E}^{(k)} = \bigcup_{ij\in E}\tilde{\sM}_{ij},$$
where $\tilde{\sM}_{ij}$ is a matching of $\{i\}\times K$ and $\{j\}\times K$.

It is easy to see that these matchings are independent, because they
correspond to successive matchings of the lifted images of distinct
edges of $G$. Moreover, two vertices $(i,\ell_1,\dots,\ell_s)\in
\{i\}\times K$ and $(j,r_1,\dots,r_s)\in j\times \{j\}\times K$ are
matched in $\tilde{\sM}_{ij}$ if  $(i,\ell_1)$ is matched to
$(j,r_1)$ in $G_1$ {\em and } $(i,r_1,r_2)$ is matched to
$(j,r_1,r_2)$ in $G_2$ {\em and} $\dots$ $(i,\ell_1,\dots,\ell_s)$
is matched to $(j,r_1,\dots,r_s)$ in $G_s$. The recipe for
constructing $G_s$ implies that the probability of this event is:
$$\Pr{\{(i,\ell_1,\dots,\ell_s),(j,r_1,\dots,r_s)\}\in\tilde{\sM}_{ij}} = \frac{1}{k_1k_2\dots k_s} = \frac{1}{k}.$$Thus if we label the elements of $K$ with the numbers $1,2,\dots,k$, we see that $\tilde{G}^{(k)}$ satisfies the assumptions of the Theorem. It follows that, just as in the case of $G^{(k)}$,
$$\Pr{\|\tilde{A}^{(k)} - A\otimes \Pi^{(k)}\|\leq 16\sqrt{\Delta\ln(4nk/\delta)}}\geq 1 - \frac{\delta}{2}.$$
Putting this together with \eqnref{closepacas} finishes the proof.\end{proof}

\section{Extensions and open questions}

{\em Lifts of Markov chains.} The argument we showed can be applied
to lifts of weighted graphs, or equivalently, of {\em reversible
Markov chains}. Let $P$ be the transition matrix of an irreducible
Markov chain on $[n]$ that is reversible with respect to a
probability measure $\pi$, meaning that $\pi(i)P(i,j)=\pi(j)P(j,i)$
for all $1\leq i,j\leq n$. [This implies that $P$ has $n$ real
eigenvalues.]

Choose a matching $\sM_{ij}$ for each pair $1\leq i\leq j\leq n$ in
the same way as in \thmref{randomlifts} and consider a Markov chain
$P^{(k)}$ on $[n]\times [k]$ with transition probabilities given by:
$$P^{(k)}((i,r),(j,\ell)) = \left\{\begin{array}{ll}P(i,j) & \{(i,r),(j,\ell)\}\in\sM_{ij};\\ 0 & \mbox{ if not}.\end{array}\right.$$
One can show (proof omitted) that the spectrum of $P^{(k)}$ contains
that of $P$ and that all new eigenvalues of $P^{(k)}$ satisfy:
$$\Pr{\max_{\eta\in {\rm new}(P^{(k)})}|\eta|\leq
16\sqrt{c_P\ln(nk/\delta)}}\geq 1-\delta,$$ where
$$c_P\equiv \max_{i\in [n]}\sum_{j=1}^n
\frac{\pi(j)P(j,i)^2}{\pi(i)}.$$ To prove this, one only needs to
consider the symmetric matrix $Q$ with entries equal to
$$Q(i,j)\equiv \sqrt{\frac{\pi(i)}{\pi(j)}}\,P(i,j)$$ (which has the same spectrum as $P$)
and the corresponding matrix $Q^{(k)}$ for the lifted chain
$P^{(k)}$, which is reversible with respect to the probability
distribution: $$\pi^{(k)}(i,\ell)=\pi(i)/k\;\;((i,\ell)\in
[n]\times[k]).$$ Notice that the parameter $c_P$ always satisfies:
$$c_P\leq \max_{i\in [n]}\left\{\,(\max_{r\in [n]}P(r,i))\sum_{j=1}^n
\frac{\pi(j)P(j,i)}{\pi(i)}\right\}= \max_{(i,r)\in [n]^2}P(r,i).$$\\

{\em Sharpness of the bound:} We do now know if the bound in
\thmref{randomlifts} can be improved. For instance, could it be the
case that all new eigenvalues of the adjacency matrix are
$\bigoh{\sqrt{\Delta}}$ with high probability, at least when the
minimum degree is $\ohmega{\ln n}$? This would be similar to the
Erd\"{o}s-R\'{e}nyi random graph \cite{FeigeOfek_Spectral} and also
related to results on random regular graphs
\cite{Friedman_SecondEigenvalue}. An analysis of the proof of
\thmref{randomlifts} shows that the only obstacle to obtaining such
a bound is the $d$ term in \thmref{freedman}, but that term is known
to be necessary in general \cite{Eu_Freedman}. However, it might be
possible to obtain better concentration bounds in the
graph-theoretic setting, at least for ``well-behaved" base graphs
$G$.

\end{document}